\documentclass[11pt]{article}
\usepackage{amssymb,amsfonts,amsmath,latexsym,epsf,tikz,url}

\newtheorem{theorem}{Theorem}[section]

\newtheorem{corollary}[theorem]{Corollary}
\newtheorem{lemma}[theorem]{Lemma}

\newcommand{\proof}{\noindent{\bf Proof.\ }}
\newcommand{\qed}{\hfill $\square$\medskip}

\begin{document}

\title{Distinguishing number and distinguishing index of neighbourhood corona  of two  graphs}

\author{
Saeid Alikhani$^{}$\footnote{Corresponding author}  \and Samaneh Soltani 
}

\date{\today}

\maketitle

\begin{center}
Department of Mathematics, Yazd University, 89195-741, Yazd, Iran\\
{\tt alikhani@yazd.ac.ir, s.soltani1979@gmail.com}
\end{center}

\begin{abstract}
	The distinguishing number (index) $D(G)$ ($D'(G)$) of a graph $G$ is the least integer $d$ such that $G$ has an vertex labeling (edge labeling)  with $d$ labels  that is preserved only by a trivial automorphism. The neighbourhood corona of two graphs $G_1$ and $G_2$ is denoted by  $G_1 \star G_2$  and is the graph obtained by
	taking one copy of $G_1$ and  $|V(G_1)|$ copies of $G_2$, and joining the neighbours of the $i$th vertex of $G_1$ to every vertex in the $i$th copy of $G_2$. In this paper we describe the automorphisms of the graph $G_1\star G_2$. Using results on automorphisms, we study the distinguishing number and the distinguishing index of $G_1\star G_2$.  We obtain upper bounds for $D(G_1\star G_2)$ and $D'(G_1\star G_2)$.  

\end{abstract}


\medskip
\noindent{\bf AMS Subj.\ Class.:} 05C15, 05E18

\section{Introduction}

Let $G = (V ,E)$ be a simple  graph with $n$ vertices.   Throughout this paper we consider only simple graphs.
The set of all automorphisms of $G$, with the operation of composition of permutations, is a permutation group
on $V$ and is denoted by $Aut(G)$.  
A labeling of $G$, $\phi : V \rightarrow \{1, 2, \ldots , r\}$, is  $r$-distinguishing, 
if no non-trivial  automorphism of $G$ preserves all of the vertex labels.
In other words,  $\phi$ is $r$-distinguishing if for every non-trivial $\sigma \in Aut(G)$, there
exists $x$ in $V$ such that $\phi(x) \neq \phi(x\sigma)$. 
The distinguishing number of a graph $G$ has defined by Albertson and Collins \cite{Albert} and  is the minimum number $r$ such that $G$ has a labeling that is $r$-distinguishing.  
Similar to this definition, Kalinkowski and Pil\'sniak \cite{Kali1} have defined the distinguishing index $D'(G)$ of $G$ which is  the least integer $d$
such that $G$ has an edge colouring   with $d$ colours that is preserved only by a trivial
automorphism.  These indices  has developed  and number of papers published on this subject (see, for example \cite{soltani2,Klavzar,fish}).

\begin{figure}[ht]
	\begin{center}
		\includegraphics[width=0.37\textwidth]{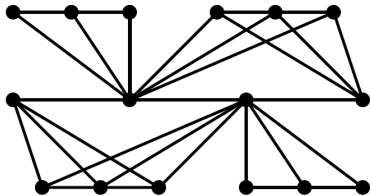}
		\caption{\label{fig1} The neighbourhood corona of $P_4 \star P_3$.}
	\end{center}
\end{figure}

We use the following notations: The set of vertices adjacent in $G$ to a vertex of a
vertex subset $W\subseteq  V$ is the open neighborhood $N_G(W )$ of $W$. The closed neighborhood  $G[W ]$ also includes all vertices of $W$ itself. In case of a singleton set $W =\{v\}$ we write $N_G(v)$ and $N_G[v]$ instead of $N_G(\{v\})$ and $N_G[\{v\}]$, respectively.  We omit the subscript when the graph $G$ is clear from the context.  The complement of  $N[v]$ in $V(G)$ is  denoted by $\overline{N[v]}$.  We denote the degree of a vertex $v$ in graph $G$ by $d_G(v)$ and the distance between two vertices $u$ and $w$ in graph $G$, by $dist_G(u,w)$. 
The corona of two graphs $G$ and $H$ which denoted by $G\circ H$ is defined in \cite{Harary} and there have been some results on the
corona of two graphs \cite{Frucht}. In \cite{soltani2} we have studied the distinguishing number and the distinguishing index of corona of two graphs. In this paper we consider another variation of corona of two graphs and study its  distinguishing number and distinguishing  index. 
Given simple graphs  $G_1$ and $G_2$, the neighbourhood corona of $G_1$ and $G_2$, denoted by $G_1 \star G_2$ and  is the graph obtained by
taking one copy of $G_1$ and $|V(G_1)|$ copies of $G_2$ and  joining the neighbours of the $i$th vertex of $G_1$ to every vertex in the $i$th copy of $G_2$ (\cite{Gopalapillai}). 
  Figure \ref{fig1} shows $P_4\star P_3$, where $P_n$ is  the path of order $n$. 
Liu and Zhu in \cite{linear} determined the adjacency spectrum of $G_1\star G_2$ for arbitrary $G_1$ and $G_2$ and the Laplacian spectrum and signless Laplacian spectrum of $G_1\star G_2$ for regular $G_1$ and arbitrary $G_2$, in terms of the corresponding spectrum of $G_1$ and $G_2$. Also Gopalapillai in \cite{Gopalapillai} has studied 
the eigenvalues and spectrum of $G_1\star G_2$, when $G_2$ is regular.

\medskip 

 In this paper we consider the neighbourhood corona  of two graphs and discuss their  distinguishing number and index.  In the next section, we give
a complete description of the automorphisms  of neighbourhood corona of two arbitrary graphs. In Section 3,  we study the distinguishing number and the  distinguishing index of neighbourhood corona of two graphs.



\section{Description of automorphisms of $G_1\star G_2$}

In this section we consider the neighbourhood corona of two graphs and describe  its automorphisms.  Let  $G_i$ has order $n_i$ and size $m_i$ $(i=1,2)$. 
The neighbourhood corona $G_1 \star G_2$ of $G_1$ and $G_2$ has $n_1+n_1n_2$ vertices and $m_1(2n_2+1) + n_1m_2$ edges and when $G_2 = K_1$, the graph $G_1 \star G_2$ is the splitting graph which has defined in \cite{Sampathkumar}.

Let $V (G_1) = \{v_1, v_2, \ldots , v_{n_1}\}$ and $V (G_2) =
\{u_1, u_2, \ldots , u_{n_2}\}$. For $i = 1, 2, \ldots , n_1$, let $u^i_1, u^i_2,\ldots , u^i_{n_2}$ denote the vertices of the $i$th copy of
$G_2$, with the understanding that $u^i_j$ is the copy of $u_j$ for each $j$. It is clear that the degrees of the vertices of $G_1 \star G_2$ are:
\begin{equation}\label{e1}
d_{G_1\star G_2}(v_i) = (n_2 + 1)d_{G_1}(v_i),~~ i = 1, 2,\ldots , n_1. 
\end{equation}
\begin{equation}\label{e2}
d_{G_1\star G_2}(u^i_j) = d_{G_2}(u_j) + d_{G_1}(v_i),~~ i = 1, 2,\ldots , n_1,~ j = 1, 2,\ldots , n_2.
\end{equation}

Now we want to know how an automorphism of $G_1 \star G_2$ acts on the vertices $G_1$ and the vertices of copies $G_2$. First we state and prove the following lemma. 
\begin{lemma}\label{lem1}
Let $G_1$ and $G_2$ be two connected graphs such that $G_1\neq K_1$ and $f$ be an automorphism of $G_1 \star G_2$ such that $f(v_i)=u^k_j$ for some $ i ,k = 1, 2,\ldots , n_1$ and $ j = 1, 2,\ldots , n_2$. Then $d_{G_1}(v_k) > d_{G_1}(v_i)$.
\end{lemma}
\proof   Since  $f(v_i)=u^k_j$, so $d_{G_1 \star G_2}(v_i)=d_{G_1 \star G_2}(u^k_j)$. By Equations (\ref{e1}) and (\ref{e2}) we have $(n_2 + 1)d_{G_1}(v_i)=d_{G_2}(u_j) + d_{G_1}(v_k)$. By contradiction, suppose  that $d_{G_1}(v_k) \leqslant d_{G_1}(v_i)$. Hence $(n_2 + 1)d_{G_1}(v_i)\leqslant d_{G_2}(u_j) + d_{G_1}(v_i)$, and so $n_2 d_{G_1}(v_i)\leqslant d_{G_2}(u_j)$.  This contradiction forces us to conclude that $d_{G_1}(v_k) > d_{G_1}(v_i)$.\qed

By Lemma \ref{lem1} we can prove the following corollary:
\begin{corollary}
Let $G_1$ be a connected graph such that $G_1\neq K_1$ and  $f$ be an arbitrary automorphism of $G_1 \star G_2$.
\begin{itemize}
\item[(i)] If  $v$ is the vertex of $G_1$ with the maximum degree in $G_1$, then $f(v)\in G_1$.
\item[(ii)] If $G_1$ is  a regular graph, then  the restriction of $f$ to $G_1$ is an automorphism of $G_1$.
\end{itemize}
\end{corollary}

We shall obtain  some results  for the automorphisms of $G_1 \star G_2$. 

\begin{figure}[ht]
	\begin{center}
		\includegraphics[width=0.6\textwidth]{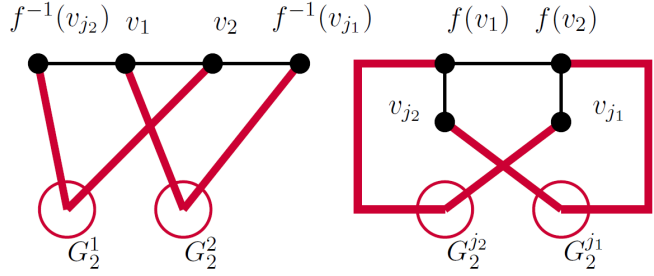}
		\caption{\label{fig14} {\small A piece of neighbourhood corona of $G_1$ and $G_2$ in the proof of Lemma \ref{lemgener1}}.}
	\end{center} 
\end{figure}

\begin{figure}[ht]
	\begin{center}
		\includegraphics[width=0.8\textwidth]{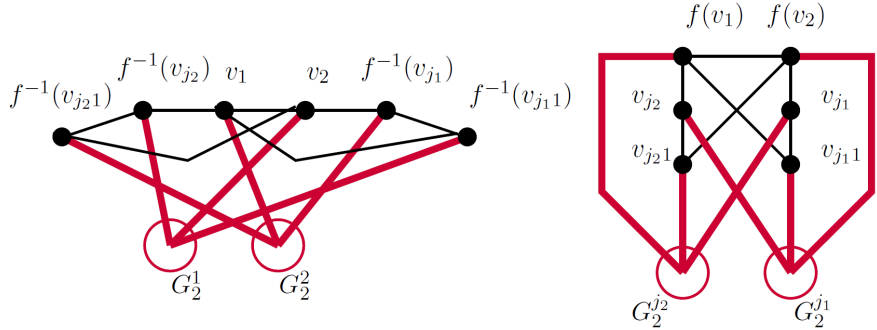}
		\caption{\label{fig15} {\small A piece of  $G_1 \star G_2$ in the proof of Lemma \ref{lemgener1}}.}
	\end{center} 
\end{figure}

\begin{lemma}\label{lemgener1}
	Let $G_1$ and $G_2$ be two connected graphs of orders $n_1$ and $n_2$, respectively, and  $n_1 > 1$. Suppose that $f$ is an automorphism of  $G_1 \star G_2$ such that the restriction of $f$ to $G_1$ is an automorphism of $G_1$, and also $f$ maps the copies of $G_2$ to each other. Then there are the automorphism $g$ of $G_1$ and the automorphisms $h_1, \ldots , h_{n_1}$ of $G_2$ such that $f(G_2^{i})=(h_i (G_2))^{k}$, where $v_k = g(v_i)$ and $i,k=1,\ldots , n_1$.
\end{lemma}
\proof
Let $f$ be an automorphism of $G_1 \star G_2$ such that the restriction of $f$ to $G_1$ is an automorphism of $G_1$, and also $f$ maps the copies of $G_2$ to each other.  Let $f$ maps the $i$th copy of $G_2$, $G_2^{i}$, to the $j_i$th copy of $G_2$, $G_2^{j_i}$, where $i,j_i=1,\ldots , n_1$, such that for the fixed numbers $i$ and $j_i$ we have $f(u_k^i)= u_{k'}^{j_i}$, where $k,k'=1,\ldots , n_2$. Then we define the automorphism $h_i$ on $G_2$ such that $h_i(u_k)=u_{k'}$. To complete the  proof we need to show that the map $g$ on $V(G_1)$ such that $g(v_i)= v_{j_i}$ is an automorphism of $G_1$, where $i,j_i=1,\ldots , n_1$. Without loss of generality we can assume that  the vertices $v_1$and $v_2$ are adjacent, and show that $v_{j_1}$ and $v_{j_2}$ are adjacent. Since the vertices $v_1$ and $v_2$ are adjacent, the vertex $v_1$ is adjacent to each vertex of $G_2^2$ (we show this concept by $v_1 \sim G_2^2$). Hence $f(v_1)\sim (h_2(G_2))^{j_2}$, and so $f(v_1)\sim v_{j_2}$ and  $v_1\sim f^{-1}(v_{j_2})$, and thus $ f^{-1}(v_{j_2})\sim G_2^1$, and finally we have $v_{j_2}\sim G_2^{j_1}$. With  a similar argument we can conclude that $f(v_2)\sim v_{j_1}$, and so $v_2\sim f^{-1}(v_{j_1})$, and hence $f^{-1}(v_{j_1})\sim G_2^2$, and thus $v_{j_1}\sim G_2^{j_2}$ (see the Figure \ref{fig14}). 
On the other hand, since $f$ maps $G_2^1$ to $(h_2(G_2))^{j_1}$, we have $d_{G_1 \star G_2}(u_k^1)=d_{G_1 \star G_2}((h_2(u_k))^{j_1})$.  We deduce  from Equations (\ref{e1}), (\ref{e2}) and $d_{G_2}(u_k)=d_{G_2}(h_2(u_k))$, that $d_{G_1}(v_1)=d_{G_1}(v_{j_1})$. Similarly, $d_{G_1}(v_2)=d_{G_1}(v_{j_2})$. Since   the restriction of $f$ to $G_1$ is an automorphism of $G_1$, we have  $d_{G_1}(v_1)=d_{G_1}(f(v_1))$ and $d_{G_1}(v_2)=d_{G_1}(f(v_2))$. Then
\begin{equation}\label{e4}
d_{G_1}(v_1)=d_{G_1}(v_{j_1})=d_{G_1}(f(v_1)),~d_{G_1}(v_2)=d_{G_1}(v_{j_2})=d_{G_1}(f(v_2)).
\end{equation}
In regard to Equation (\ref{e4}) and Figure \ref{fig14}, there exists the vertices $v_{j_11}$ and $v_{j_21}$ adjacent to vertices $v_{j_1}$ and $v_{j_2}$, respectively. Thus the vertices $v_{j_11}$ and $v_{j_21}$ are adjacent to $G_2^{j_1}$ and $G_2^{j_2}$, respectively, and so $f^{-1}(v_{j_11})\sim G_2^1$ and $f^{-1}(v_{j_21})\sim G_2^2$. Hence $f^{-1}(v_{j_11})\sim v_1$ and $f^{-1}(v_{j_21})\sim v_2$. Since $v_{j_1}\sim v_{j_11}$ and $v_{j_2}\sim v_{j_21}$, so $f^{-1}(v_{j_1})\sim f^{-1}(v_{j_11})$ and $f^{-1}(v_{j_2})\sim f^{-1}(v_{j_21})$ (see Figure \ref{fig15}).
Note that, for every vertex in $N_G(v_{j_2})$ such as $x$, we have $x\sim G_2^{j_2}$. So we see that $f^{-1}(x)\sim G_2^2$, and so $f^{-1}(x)\sim v_2$ (similar argument satisfies for each vertex in $N_G(v_{j_1})$). In regard to Figure \ref{fig15} and Equation (\ref{e4}), we need to the other vertex adjacent to $v_{j_1}$, such as $x$. If $x$ has been chosen among the nonadjacent vertices to $G_2^{j_2}$ that has been shown in Figure \ref{fig15}, then with the similar argument as above, we obtain that $f^{-1}(x)$ is adjacent to $v_1$, and so Equation (\ref{e4}) dose  not satisfy, again. Therefore after finite steps we should choose a vertex adjacent to $v_{j_1}$, such as $x$, among the vertices that are adjacent to $G_2^{j_2}$, otherwise we conclude that  the order of $G_1$ is infinite and this is a  contradiction.  By Figure \ref{fig15} and above  information, the vertex $v_{j_1}$ is the only vertex that is adjacent to $G_2^{j_2}$ and is not among the adjacent vertices to $v_{j_1}$, in each step. Hence $v_{j_1}\sim v_{j_2}$, and the result follows.\qed

\begin{figure}[ht]
	\begin{center}
		\includegraphics[width=0.8\textwidth]{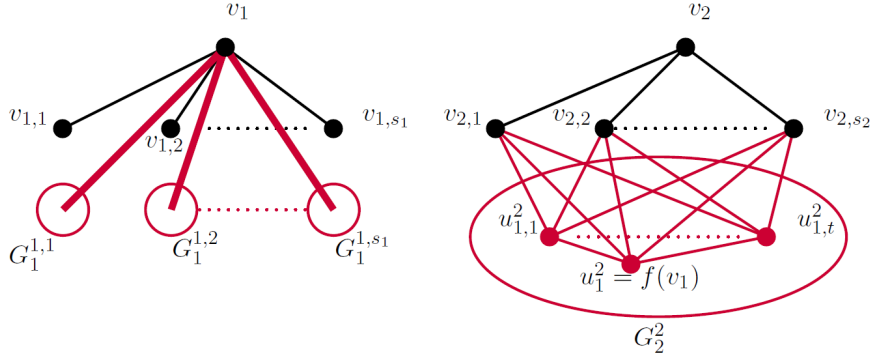}
		\caption{\label{fig16} {\small A piece of neighbourhood corona of $G_1$ and $G_2$ in the proof of Lemma \ref{lemmgener2}}.}
	\end{center} 
\end{figure}

\begin{lemma}\label{lemmgener2}
	Let $G_1$ and $G_2$ be two connected graphs of order $n_1$ and $n_2$, respectively, and  $n_1 > 1$. If $f$ is an automorphism of  $G_1 \star G_2$, then the restriction of $f$ to $G_1$ is an automorphism of $G_1$.
\end{lemma}
\proof
Since $f$ is an automorphism, it is suffices  to show that  the restriction of $f$ to $G_1$ is an automorphism of $G_1$. By contradiction, suppose   that $f\vert_{G_1}$ is not an automorphism of $G_1$. Without loss of generality we assume that $f(v_1)=u_1^2$. Hence by  Lemma \ref{lem1}, $d_{G_1}(v_2)> d_{G_1}(v_1)$. Since $f$ preserves the degree of the vertices, $d_{G_1 \star G_2}(v_1)= d_{G_1 \star G_2}(u_1^2)$, and so by Equations (\ref{e1}) and (\ref{e2}) we have $(1+n_2)d_{G_1}(v_1)=d_{G_2}(u_1)+d_{G_1}(v_2)$. Suppose that $N_{G_1}(v_1)=\{v_{1,1},\ldots, v_{1,s_1}\}$, $N_{G_1}(v_2)=\{v_{2,1},\ldots, v_{2,s_2}\}$ and $N_{G_2}(u_1)=\{u_{1,1},\ldots, u_{1,t}\}$ where  $(1+n_2)s_1=t+s_2$ and  $s_i=d_{G_1}(v_i)$, $i=1,2$, and also $t=d_{G_2}(u_1)$ (see Figure \ref{fig16}). 
Since  $f$ preserves the adjacency relation, so $f(N_{G_1 \star G_2}(v_1))=N_{G_1 \star G_2}(u_1^2)$, i.e.,
\begin{align}\label{e5}
\big\{f(v_{1,1}),\ldots , f(v_{1,s_1}),f(u_1^{1,1}),\ldots, f(u_{n_2}^{1,1}),&\ldots , f(u_{1}^{1,s_1}),\ldots , f(u_{n_2}^{1,s_1}) \big\}\nonumber\\
&=\{u_{1,1}^2,\ldots , u_{1,t}^2,v_{2,1},\ldots , v_{2,s_2}\}.
\end{align}
Since $t< n_2$, there are  vertices in the copies $G_2^{1,1},\ldots , G_2^{1,s_1}$ such that they are mapped to the elements of the set $\{v_{2,1},\ldots , v_{2,s_2}\}$, under the automorphism $f$. Without loss of generality we can assume that $f(u_{i_j}^{1,j})=v_{2,j}$, where $1\leqslant j\leqslant s_1$. We continue the proof by considering two cases for $s_1$ as follows:

Case 1) If  $s_1> 1$. Since $v_2$ is adjacent to the vertices $v_{2,1},\ldots , v_{2,s_1}$, so $f^{-1}(v_2)$ is adjacent to the vertices $u_{i_1}^{1,1},\ldots , u_{i_{s_1}}^{1,s_1}$. Since $s_1> 1$, so $f^{-1}(v_2)\in G_1$ and $f^{-1}(v_2)$ is adjacent to the vertices $v_{1,1},\ldots , v_{1,s_1}$. Hence $v_2$ is adjacent to the vertices $f(v_{1,1}),\ldots , f(v_{1,s_1})$, and by Equation \ref{e5} we have 
\begin{equation}\label{e6}
\{f(v_{1,1}),\ldots , f(v_{1,s_1})\}\subseteq \{v_{2,s_1+1},\ldots , v_{2,s_2}\}.
\end{equation}

\begin{figure}[ht]
	\begin{center}
		\includegraphics[width=1.0\textwidth]{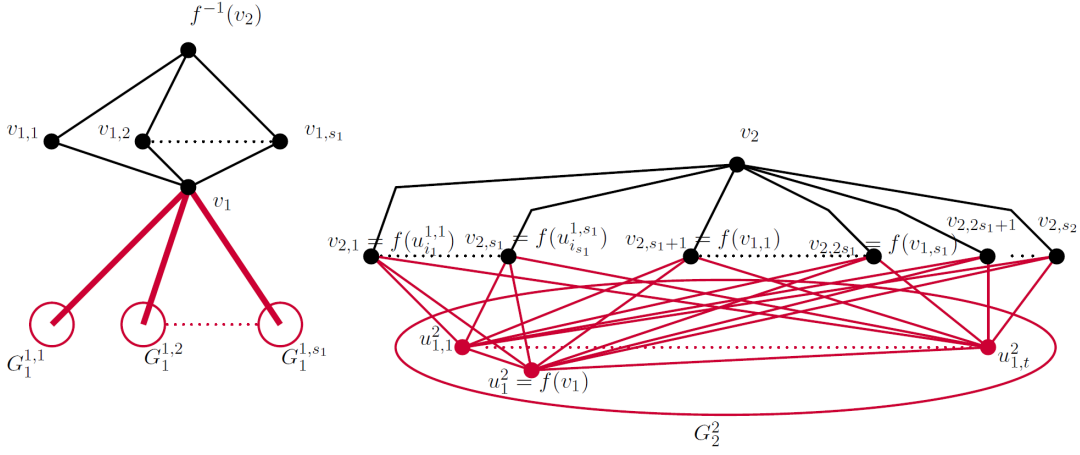}
		\caption{\label{fig17} {\small A piece of neighbourhood corona of $G_1$ and $G_2$ in the proof of Lemma \ref{lemmgener2}}.}
	\end{center} 
\end{figure}

Without loss of generality we  assume  that $f(v_{1,i})=v_{2,s_1+i}$, where $1 \leqslant i \leqslant s_1$ (see Figure \ref{fig17}).

Since $f^{-1}(v_2)$ is adjacent to the vertices $v_{1,1},\ldots , v_{1,s_1}$, we can say that $f^{-1}(v_2)$ is adjacent to all vertices of $G_2^{1,1},\ldots , G_2^{1,s_1}$, so $v_2$ is adjacent to all vertices of $f(G_2^{1,1}),\ldots , f(G_2^{1,s_1})$. Then by Equation (\ref{e5}) we get 
\begin{equation}\label{e7}
\{f(u_1^{1,1}),\ldots, f(u_{n_2}^{1,1}),\ldots , f(u_{1}^{1,s_1}),\ldots , f(u_{n_2}^{1,s_1})\}\subseteq \{v_{2,1},\ldots , v_{2,s_2}\}.
\end{equation}

With respect to Equations (\ref{e5}), (\ref{e6}) and (\ref{e7}) we have a contradiction.

\medskip

\begin{figure}[ht]
	\begin{center}
		\includegraphics[width=0.8\textwidth]{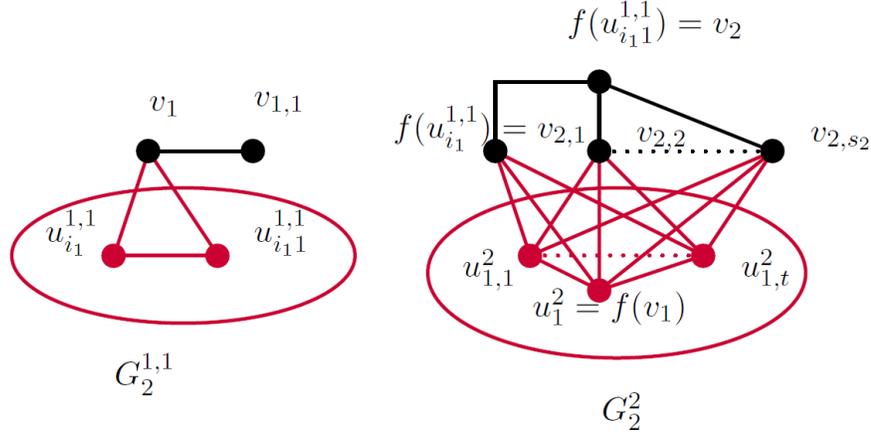}
		\caption{\label{fig18} {\small A piece of neighbourhood corona of $G_1$ and $G_2$ in the proof of Lemma \ref{lemmgener2}}.}
	\end{center} 
\end{figure}

Case 2) If  $s_1=1$. Since $f$ preserves the adjacency relation,  so 
\begin{equation}\label{e8}
\{f(v_{1,1}), f(u_1^{1,1}),\ldots, f(u_{n_2}^{1,1})\}= \{u_{1,1}^2,\ldots, u_{1,t}^2, v_{2,1},\ldots , v_{2,s_2}\}.
\end{equation} 
Since  $t< n_2$, there exists a vertex in the copy $G_2^{1,1}$ such that it is mapped to an elements of the set $\{v_{2,1},\ldots , v_{2,s_2}\}$, under the automorphism $f$. Without loss of generality we can assume that $f(u_{i_1}^{1,1})=v_{2,1}$. Since $v_2$ is adjacent to $v_{2,1}$, so $f^{-1}(v_2)$ is adjacent to $u_{i_1}^{1,1}$, and since $f^{-1}(v_2)\neq v_1$, so $f^{-1}(v_2)\in G_2^{1,1}$. Without loss of generality we can assume that $f^{-1}(v_2)=u_{i_11}^{1,1}$ such that $u_{i_11}^{1,1}$ is adjacent to $u_{i_1}^{1,1}$ (see Figure \ref{fig18}).

Since $v_{1,1}$ is adjacent to the vertex $v_1$ and $dist_{G_1 \star G_2}(v_{1,1},u_{i_1}^{1,1})=dist_{G_1 \star G_2}(v_{1,1},u_{i_11}^{1,1})=2$, so $f(v_{1,1})$ is adjacent to the vertex $u_1^2$ and also
\begin{equation}\label{e9}
dist_{G_1 \star G_2}(f(v_{1,1}),v_2)=dist_{G_1 \star G_2}(f(v_{1,1}),v_{2,1})=2.
\end{equation}

Now by Equations (\ref{e8}) and (\ref{e9}) we have a contradiction. Therefore the restriction of each automorphism of $G_1 \star G_2$ to $G_1$ is an automorphism of $G_1$.\qed

\begin{corollary}\label{cor}
	Let $G_1$ and $G_2$ be two connected graphs of order $n_1$ and $n_2$, respectively, such that $n_1 > 1$ and  $f$ is an automorphism of  $G_1 \star G_2$. Then the restriction of $f$ to $G_1$ is an automorphism of $G_1$ and also there are the automorphism $g$ of $G_1$ and the automorphisms $h_1, \ldots , h_{n_1}$ of $G_2$ such that $f(G_2^{i})=(h_i (G_2))^{k}$, where $v_k = g(v_i)$ and $i,k=1,\ldots , n_1$.
\end{corollary}
\proof
By Lemmas  \ref{lemgener1} and \ref{lemmgener2}, it is sufficient to prove that the copies of $G_2$ are mapped to each other under the automorphism $f$, and it follows from that  $f$ preserves the adjacency relation on each copy of $G_2$. \qed

The following corollary is an immediate consequence of Corollary \ref{cor} for 
graphs of the form $G\star K_1$. 

\begin{corollary}\label{lem2}
	Let $G$ be a connected graph of order $n>1$ and  $f$ be an arbitrary automorphism of $G \star K_1$. Then the restriction of $f$ to $G$ is an automorphism of $G$. Also  $f(K_1^i)=K_1^{j_i}$ for some automorphism $g$  of $G$ such that $g(v_i)=v_{j_i}$ where $ i ,j_i = 1, 2,\ldots , n_1$.
\end{corollary}

\section{Study of $D(G_1\star G_2)$ and $D'(G_1\star G_2)$}

In this section we use the results in Section 2 to study the distinguishing number and the distinguishing index of the neighbourhood corona of two graphs. First we consider the neighbourhood corona of an arbitrary graph with $K_1$. 
The following theorem gives an upper bound for $D(G\star K_1)$ and $D'(G\star K_1)$.  
\begin{theorem}\label{displitt}
Let $G$ be a connected graph of order $n>1$. We have 
\begin{itemize}
\item[(i)] $D(G\star K_1)\leqslant D(G)$,
\item[(ii)] $D'(G\star K_1)\leqslant D'(G)$.
\end{itemize}
\end{theorem}
\proof
\begin{enumerate} 
	\item[(i)]  We shall define a distinguishing vertex labeling for $G\star K_1$ with $D(G)$ labels. First we label $G$ in a distinguishing way with $D(G)$ labels. Next we assign the vertex $K_1^{v_i}$, the label of the vertex $v_i$ where $1\leqslant i \leqslant n$. This labeling is a distinguishing vertex labeling of $G\star K_1$, because if $f$ is an automorphism of $G\star K_1$ preserving the labeling then by Corollary \ref{cor}, the restriction of $f$ to $G$ is an automorphism of $G$ preserving the labeling. Since we labeled $G$ in a distinguishing way at first, so  the restriction of $f$ to $G$ is the identity automorphism on $G$. On the other hand by Corollary \ref{lem2} there exists an automorphism $g$ of $G$ such that $f(K_1^{v_i})=K_1^{g(v_i)}$, $1\leqslant i \leqslant n$. Regarding to the  labeling of copies of $K_1$, we can obtain that $g$ is the identity automorphism on $G$, and so  $f$ is the identity automorphism on $G\star K_1$.

\item[(ii)] 
 We define a distinguishing edge labeling for $G\star K_1$ with $D'(G)$ labels. First we label the edges of $G$ in a distinguishing way with $D'(G)$ labels. By Equations (\ref{e1}) and (\ref{e2}) we know that the degree of $K_1^{v_i}$ in $G\star K_1$ is equal with the degree of $v_i$ in $G$ where $1\leqslant i \leqslant n$.  Now we assign the edge between  $K_1^{v_i}$ and $v_{i,j}$ where $v_{i,j}\in N_G(v_i)$, the label of the edges between $v_i$ and $v_{i,j}$ where $j=1,\ldots , d_G(v_i)$.   This labeling is a distinguishing edge labeling of $G\star K_1$, because if $f$ is an automorphism of $G\star K_1$ preserving the labeling then by Corollary \ref{cor}, the restriction of $f$ to $G$ is an automorphism of $G$ preserving the labeling. Since we labeled $G$ in a distinguishing way at first, so  the restriction of $f$ to $G$ is the identity automorphism on $G$. On the other hand by Corollary \ref{lem2} there exists an automorphism $g$ of $G$ such that $f(K_1^{v_i})=K_1^{g(v_i)}$, $1\leqslant i \leqslant n$. Regarding to the  labeling of the edges incident to each copies of $K_1$, we can obtain that $g$ is the identity automorphism on $G$, and so  $f$ is the identity automorphism on $G\star K_1$.\qed
\end{enumerate}

\begin{figure}[ht]
\begin{center}
\includegraphics[width=0.45\textwidth]{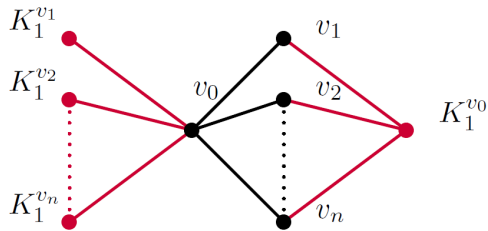}
	\caption{\label{fig8} {\small The  neighbourhood corona of $K_{1,n}$ and $K_1$.}}
\end{center} 
\end{figure}

The bounds of $D(G\star K_1)$ and $D'(G\star K_1)$ in Theorem \ref{displitt} are sharp. If we consider $G$ as the star graph $K_{1,n}$, $n>1$, then $K_{1,n}\star K_1$ is a graph as shown in Figure \ref{fig8}. Using  the degree of the verices of  $K_{1,n}\star K_1$ we can get the automorphism group of $K_{1,n}\star K_1$ and then it can be concluded that $D(K_{1,n}\star K_1)=n=D(K_{1,n})$, and also $D'(K_{1,n}\star K_1)=n=D'(K_{1,n})$.

\medskip 
In Theorem \ref{displitt}, the  sharp upper bounds for $D(G\star K_1)$ and $D'(G\star K_1)$ have been given, but we did not present  lower bounds for these parameters. Actually, there are graphs whose  distinguishing number can be arbitrarily larger than the   distinguishing number of its neighbourhood corona with $K_1$. In other words, we can  show that there exists a  connected graph $G$ of order $n>1$  such that the value of $\frac{D(G\star K_1)}{D(G)}$ can be arbitrarily small. To do this we
need  the two following theorems. Recall that the friendship graph $F_n$ $(n\geqslant 2)$ can be constructed by joining $n$ copies of the cycle graph $C_3$ with a common vertex. 
 
 \begin{theorem}{\rm \cite{soltani2}}\label{disfan}
The distinguishing number of the friendship graph $F_n$  $(n\geq 2)$ is  $$D(F_n)= \big\lceil \frac{1+\sqrt{8n+1}}{2}\big\rceil .$$  
\end{theorem}

\begin{figure}[ht]
\begin{center}
\includegraphics[width=0.75\textwidth]{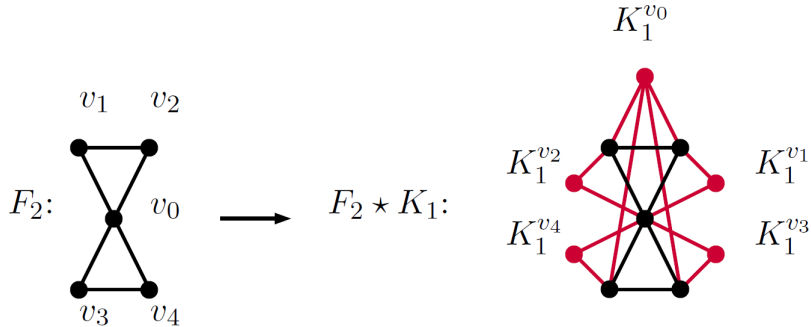}
	\caption{\label{fig9} {\small The graphs $F_2$ and $F_2 \star K_1$.}}
\end{center} 
\end{figure}

Now we obtain the exact value of the distinguishing number of neighborhood corona of $F_n$ with $K_1$. 
\begin{theorem}\label{disfanstar}
The distinguishing number of $F_n\star K_1$  $(n\geq 2)$ is  $$D(F_n\star K_1)= \big\lceil \sqrt{ \dfrac{1+\sqrt{8n+1}}{2}}\big\rceil .$$  
\end{theorem}
\proof
Let $V(F_n)=\{v_0,v_1,v_2,\ldots , v_{2n-1},v_{2n}\}$ and the vertex $v_0$ be the central vertex and $v_{2i-1}$ and $v_{2i}$ be the vertices of the base of triangles in $F_n$ where $1\leqslant i \leqslant n$. So $d_{F_n}(v_0)=2n$ and $d_{F_n}(v_i)=2$ where $1\leqslant i \leqslant 2n$. By Equations (\ref{e1}) and (\ref{e2}) we have $d_{F_n\star K_1}(v_0)=4n$ and $d_{F_n\star K_1}(v_i)=4$, also $d_{F_n\star K_1}(K_1^{v_0})=2n$ and $d_{F_n\star K_1}(K_1^{v_i})=2$ where $1\leqslant i \leqslant 2n$ (see Figure \ref{fig9}).

If $f$ is an automorphism of $F_n \star K_1$, then $f$ fixes the vertices $v_0$ and $K_1^{v_0}$ (if $n=2$ we can get the same result by Corollary \ref{cor}). So we assign the vertices $v_0$ and $K_1^{v_0}$ the label $1$. Let $(x_i,y_i,z_i,w_i)$ be the label of the vertices $(K_1^{v_{2i}}, v_{2i-1}, v_{2i}, K_1^{v_{2i-1}})$ where $1\leqslant i \leqslant n$. Suppose that $L=\{(x_i,y_i,z_i,w_i)~\vert ~1\leqslant i \leqslant n ,  x_i,y_i,z_i,w_i\in \mathbb{N}\}$, is a labeling of the vertices of $F_n \star K_1$ except the vertices $v_0$ and $K_1^{v_0}$. If $L$ is a distinguishing labeling of $F_n \star K_1$ then:
\begin{itemize}
\item[(i)] For every $i=1,\ldots, n$, it should be satisfied that $x_i\neq w_i$ or $y_i\neq z_i$. Otherwise, the automorphism $f_i$ of $F_n \star K_1$ such that $f_i$ maps $K_1^{v_{2i}}$ and $K_1^{v_{2i-1}}$ to each other, two vertices  $v_{2i-1}$ and $v_{2i}$ to each other, and fixes the remaining vertices, preserves the labeling.
\item[(ii)] For every $i$ and $j$ in $\{1,\ldots, n\}$, with $i\neq j$,  it should be satisfied that $(x_i,y_i,z_i,w_i)\neq (x_j,y_j,z_j,w_j)$ and $(x_i,y_i,z_i,w_i)\neq (w_j,z_j,y_j,x_j)$. Otherwise, the automorphism $f_{ij}$ and $g_{ij}$ of $F_n \star K_1$ by the following definitions preserve the labeling.
\item The automorphism $f_{ij}$ maps $K_1^{v_{2i}}$ and $K_1^{v_{2j}}$ to each other and also $K_1^{v_{2i-1}}$ and $K_1^{v_{2j-1}}$ to each other. The map $f_{ij}$ maps $v_{2i}$ and $v_{2j}$ to each other, also it maps $v_{2i-1}$ and $v_{2j-1}$ to each other and fixes the remaining vertices of $F_n \star K_1$.
\item The automorphism $g_{ij}$ maps $K_1^{v_{2i}}$ and $K_1^{v_{2j-1}}$ to each other, also $K_1^{v_{2i-1}}$ and $K_1^{v_{2j}}$ to each other. The map $g_{ij}$ maps $v_{2i}$ and $v_{2j-1}$ to each other, also it maps $v_{2i-1}$ and $v_{2j}$ to each other and fixes the remaining vertices of $F_n \star K_1$.
\end{itemize}

\begin{figure}[ht]
\begin{center}
\includegraphics[width=0.5\textwidth]{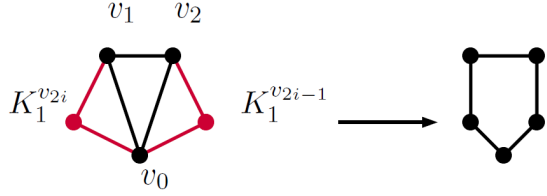}
	\caption{\label{fig10} {\small The considered  pentagon (or a cycle of size $5$) in the proof of Theorem \ref{disfanstar}}.}
\end{center} 
\end{figure}

\begin{figure}[h]
	\begin{minipage}{6.7cm}
		\includegraphics[width=\textwidth]{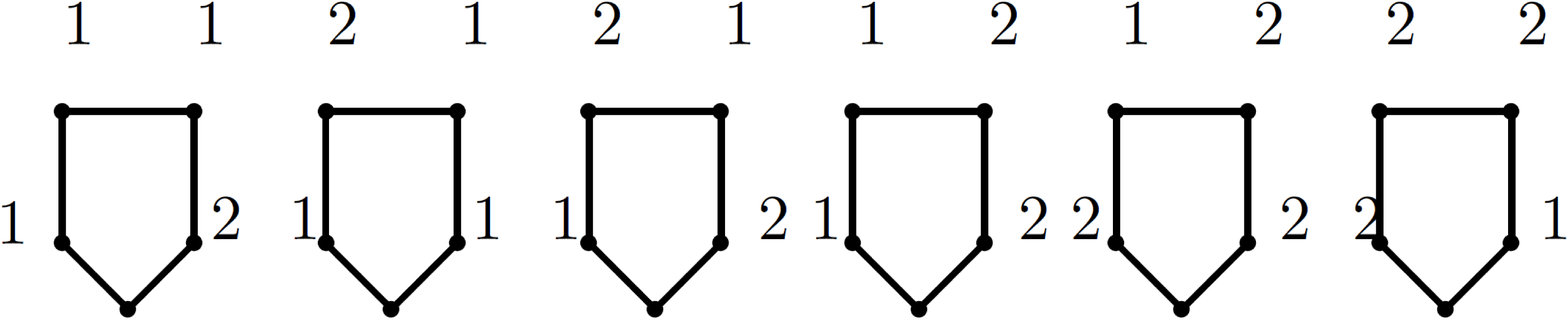}
	\end{minipage}
	\hspace{.5cm}
	\begin{minipage}{6.7cm}
		\includegraphics[width=\textwidth]{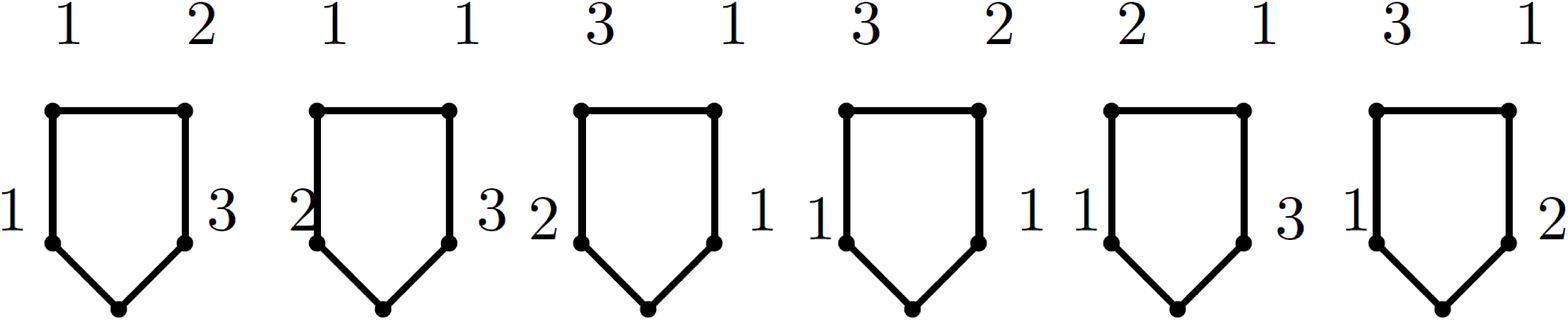}
	\end{minipage}
	\caption{\label{fig12} {\small Distinguishing labeling of blades  with the labels $\{1,2\}$ and $\{1,2,3\}$}, respectively.}
\end{figure}

So using the label set $\{1,\ldots , s\}$ we can make at most $(s^4-s^2)/2$ of the $4$-ary's $(x,y,z,w)$ satisfying $(i)$ and $(ii)$. Because, the number of $4$-ary's $(x,y,z,w)$ such that $x\neq w$ is $s(s-1)s^2$, and the number of $4$-ary's $(x,y,z,w)$ such that $y\neq z$ is $s(s-1)s^2$. On the other hand the number of $4$-ary's $(x,y,z,w)$ such that $x\neq w$ and $y\neq z$ is $(s(s-1))^2$. So the maximum number of $4$-ary's $(x,y,z,w)$ satisfying $(i)$ is $$(s(s-1)s^2 + s(s-1)s^2)-(s(s-1))^2=s^4-s^2.$$
 Among these  $4$-ary's we should choose the $4$-ary's that  satisfying $(ii)$, too. Therefore the number of $4$-ary's $(x,y,z,w)$ satisfying $(i)$ and $(ii)$ which they can make by the label set $\{1,\ldots , s\}$  is $\frac{s^4-s^2}{2}$. Therefore  $D(F_n \star K_1)\geqslant min\{s: \frac{s^4-s^2}{2} \geqslant n\}$. By an easy computation, we see that $$min\{s: \frac{s^4-s^2}{2} \geqslant n\}=\lceil \sqrt{ \dfrac{1+\sqrt{8n+1}}{2}}\rceil.$$  Now we present a distinguishing vertex labeling with this number of labels.
We assign $v_0$ and $K_1^{v_0}$ the label $1$. We should label the remaining vertices such that the identity automorphism preserves the labeling only. Denoting each pentagon with the vertices $K_1^{v_{2i}}, v_{2i-1}, v_{2i}, K_1^{v_{2i-1}},v_0$ in $F_n \star K_1$ where $1\leqslant i \leqslant n$,  by a general pentagon that have  shown in Figure \ref{fig10} and calling it a blade and  continue the labeling. 
At first, we want to know the maximum number of blades that can be labeled in a distinguishing way by $1$ and $2$. As we can see in  Figure \ref{fig12}, the maximum number of blades that can be  labeled in distinguishing way,  by $1,2$ is $\underline{6}$.

 In order to preserve the labeling under the identity automorphism only, we should use another label to assign the next blade. As mentioned earlier, the maximum number of blades that can be labeled by each the set $\{1,3\},\{2,3\}$ is six. Now we want to know the maximum number of blades that can be labeled by presence of $\{1,2,3\}$ at the same time in the blade. This number is $18$. Because let to label  with the labels $1,2,3$ and a repetition of $1$. As shown in Figure \ref{fig12},  we can label six blades. Obviously we can do the same with letting repetition of $2$ and $3$.   
Therefore the maximum number of blades that can be labeled by presence of $\{1,2,3\}$ at the same time is $\underline{18}$. Until now, we labeled $36$ blades.
\begin{equation*}
\underbrace{6}_{\{1,2\}}+\underbrace{6}_{\{1,3\}}+\underbrace{6}_{\{2,3\}}+\underbrace{18}_{\{1,2,3\}}=36
\end{equation*}

\begin{figure}[ht]
\begin{center}
\includegraphics[width=0.6\textwidth]{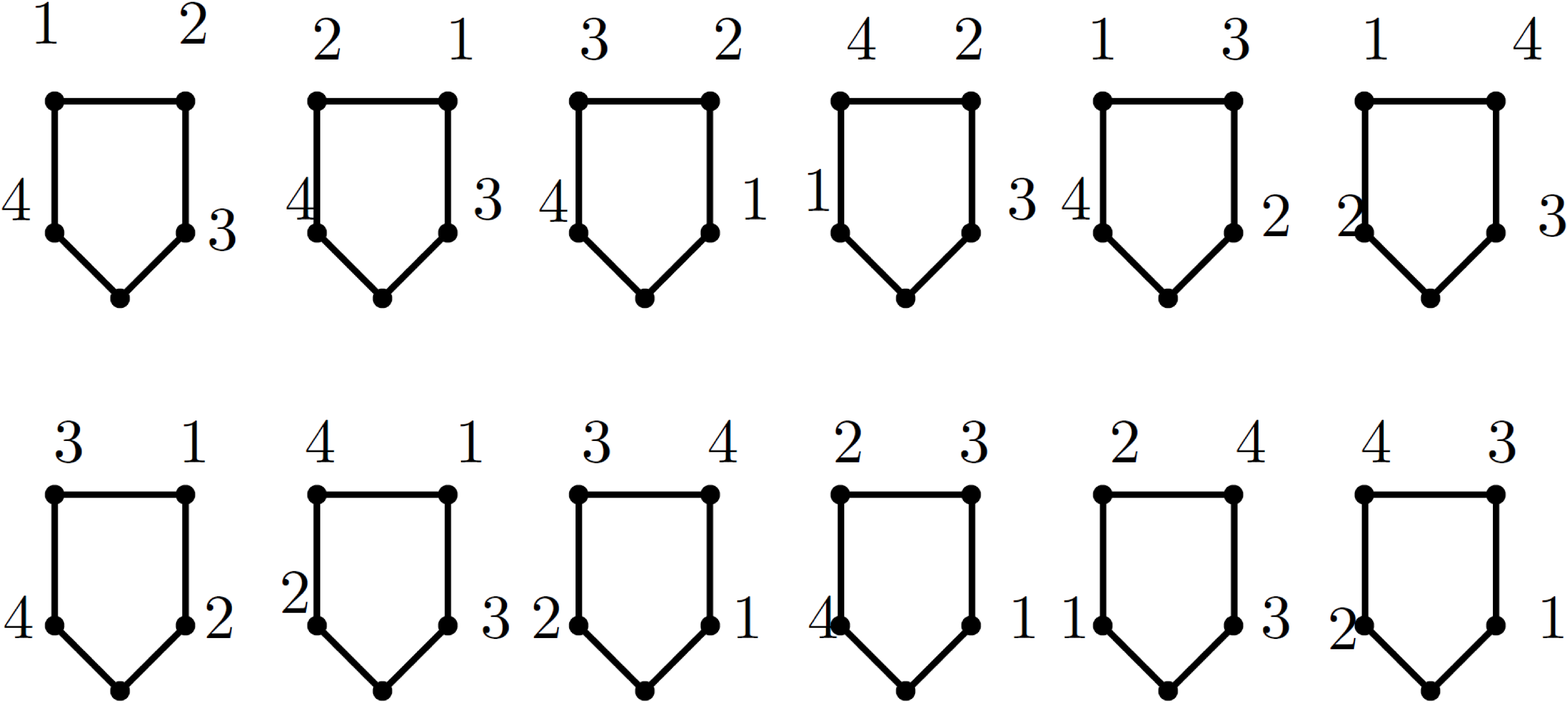}
	\caption{\label{fig13} {\small The distinguishing labeling of blades  with the labels $\{1,2,3,4\}$}.}
\end{center} 
\end{figure}

If we want to label the next blade, we should add a new label, $4$. The maximum number of blades that can be labeled by each the set $\{1,4\},\{2,4\}, \{3,4\}$ is six. Also, the maximum number of blades that can be labeled by each the set $\{1,2,4\},\{1,3,4\}, \{2,3,4\}$ is eighteen. We can see that  the maximum number of blades that can be labeled by presence of $\{1,2,3,4\}$ at the same time is  $\underline{12}$ as Shown in Figure \ref{fig13}.

Thus  we have labeled $120$ blades until now.
\begin{equation*}
36+ \underbrace{6}_{\{1,4\}}+\underbrace{6}_{\{2,4\}}+\underbrace{6}_{ \{2,4\}}+\underbrace{18}_{\{1,2,4\}}+\underbrace{18}_{ \{1,3,4\}}+\underbrace{18}_{\{2,3,4\}} + \underbrace{12}_{\{1,2,3,4\}}=120.
\end{equation*}

Therefore the relationship between the number of labels that has been used, $ \mathbf{d}(F_n\star K_1)$, and $n$  are as the following sequence:
\begin{equation*}
\{\mathbf{d}(F_n\star K_1)\}=\{0,\underbrace{2}_{6-times}, \underbrace{3}_{30-times}, \underbrace{4}_{84-times},\ldots , m, \ldots , m, \ldots \}.
\end{equation*}
where the number of the repetitions $m$ in above sequence is $ (m-1)6+{m-1 \choose 2}18+ {m-1 \choose 3}12$, with $  m\geqslant 1$.

In fact,
$ \mathbf{d}(F_n\star K_1)= min\{k:~\sum_{i=1}^{k}\left( {i-1\choose 1}6+{m-1\choose 2}18+{m-1\choose 3}12 \right) \geqslant n \}$.  By an easy computation, we see that 
\begin{align*}
&min\{k:~\sum_{i=1}^{k}\left( {i-1\choose 1}6+{m-1\choose 2}18+{m-1\choose 3}12 \right) \geqslant n \}\\
&=min\{k: (k^4-k^2)/2 \geqslant n\}\\
&=\lceil \sqrt{ \dfrac{1+\sqrt{8n+1}}{2}}\rceil.
\end{align*}
Therefore we have the result.\qed

Now we are ready to state and prove the following theorem: 
\begin{theorem}\label{thmsmall}
	There exists a  connected graph $G$ of order $n>1$  such that the value of $\frac{D(G\star K_1)}{D(G)}$ can be arbitrarily small.
\end{theorem}
\proof 
 By Theorems \ref{disfan} and \ref{disfanstar} it can be seen that 
\begin{equation*}
lim_{n\rightarrow \infty} \dfrac{D(F_n\star K_1)}{D(F_n)}=lim_{n\rightarrow \infty}\dfrac{\lceil \sqrt{ \dfrac{1+\sqrt{8n+1}}{2}}\rceil}{\lceil \dfrac{1+\sqrt{8n+1}}{2}\rceil}=0
\end{equation*}
Therefore we have the result.\qed

 The following theorem is one of the  main result of this paper and gives an upper bound for the distinguishing number of the neighbourhood corona of two arbitrary graphs: 
 
 \begin{theorem}\label{thm}
 Let $G_1$ and $G_2$ be two connected graphs of orders $n_1$ and $n_2$, respectively, such that $n_1 > 1$. Then $ D(G_1\star G_2)\leqslant max\{D(G_1), D(G_2)+M\}$,
 where $$M= min\left \{k: \sum_{m=0}^{k}y_m  \geqslant D(G_1)\right\}, ~~~ y_m=\left\{
 \begin{array}{ll}
 1&m=0,\\
 D(G_2)& m=1,\\
 D(G_2)+\sum_{i=1}^{m-1}{m-1 \choose i}{D(G_2)\choose i+1}& m\geqslant 2.
 \end{array}\right.$$
 \end{theorem}
 \proof
 We define a distinguishing vertex labeling for $G_1\star G_2$ with $max\{D(G_1), D(G_2)+M\}$ labels. First we label $G_1$ with $D(G_1)$ labels in a distinguishing way. For the labeling of copies of $G_2$, we partition the vertices of $G_1$ by the distinguishing labeling of $G_1$, i.e., we partition the vertices of $G_1$ into $D(G_1)$ classes, such that $[i]$th class contains the vertices of $G_1$ having the label $i$, in the distinguishing labeling of $G_1$, where $1\leqslant i\leqslant D(G_1)$. Let $[i]=\{v_{i1},\ldots , v_{is_i}\}$,  where $s_i$ is the size of $[i]$th class and $1\leqslant i\leqslant D(G_1)$. By this partition we label the copies of $G_2$ as follows: First we label the vertices of $G_2$ with $D(G_2)$ labels in a distinguishing way, next we do the following changes on the labeling of $G_2$. Before the labeling of the copies of $G_2$, we introduce the notation $G_2^{[i]}$ for the set $\{G_2^{i1},\ldots , G_2^{is_i}\}$, i.e., $G_2^{[i]}$ is the set of that copies of $G_2$ corresponding to the elements of $[i]$th class, where $1\leqslant i\leqslant D(G_1)$. In fact we partition the copies of $G_2$ into $D(G_1)$ classes, that $G_2^{[i]}$ is the notation of $[i]$th class. Now we present the labeling of copies of $G_2$ by the following steps:
 
 Step 1) We label all of the copies of $G_2$ which  are in $G_2^{[1]}$, exactly the same as the distinguishing labeling of $G_2$.
 
 Step 2) For the labeling of the copies in $G_2^{[i]}$, where $2\leqslant i \leqslant D(G_2)+1$, we use of the new label $D(G_2)+1$ in such a way that the label $i-1$ in the all elements of $G_2^{[i]}$ is replaced by the new label $D(G_2)+1$,   where $2\leqslant i \leqslant D(G_2)+1$.
 
 Step 3)  For the labeling of the copies in $G_2^{[i]}$, where $D(G_2)+2\leqslant i \leqslant 2D(G_2)+1$, we do the same action  as Step 2, with the new label $D(G_2)+2$, instead of the labels $D(G_2)+1$.
 
 Step 4) By choosing two labels among the labels $\{1,\ldots , D(G_2)\}$, and replacing them by the two new labels $D(G_2)+1$ and $D(G_2)+2$, we can label the elements of ${D(G_2) \choose 2}$ other classes of the classes  $G_2^{[i]}$.
 
 Step 5) We do the same work as Step 2 with the new label $D(G_2)+3$ instead of labels $D(G_2)+1$. Next we label $2{D(G_2) \choose 2}$ other classes   $G_2^{[i]}$, with the two new labels $D(G_2)+1$ and $D(G_2)+3$, also with the labels $D(G_2)+2$ and $D(G_2)+3$, exactly the same as Step 4.
 
 Step 6) Now we choose three labels among the labels $\{1,\ldots , D(G_2)\}$, and replace them by the three new labels $D(G_2)+1$, $D(G_2)+2$ and $D(G_2)+3$.
 
 By continuing this method we conclude  that the number of classes can be labeled with the labels $1,\ldots , D(G_2)+m$, $m\geqslant 1$, such that the label $D(G_2)+m$ is used in the labeling of each element of classes, is $y_m$ where
 
 \begin{equation*}
 y_m=\left\{
 \begin{array}{ll}
 1&m=0,\\
 D(G_2)& m=1,\\
 D(G_2)+\sum_{i=1}^{m-1}{m-1 \choose i}{D(G_2)\choose i+1}& m\geqslant 2.
 \end{array}\right.
 \end{equation*}
  
  Therefore the number of labels that have been used for the labeling of all copies of $G_2$, is $D(G_2)+M$ where $M= min\left \{k: \sum_{m=0}^{k}y_m  \geqslant D(G_1)\right\}$. 
   This labeling is a distinguishing vertex labeling of $G_1\star G_2$, because if $f$ is an automorphism of $G_1\star G_2$ preserving the labeling, then by Corollary \ref{cor}, $f\vert_{G_1}$ is an automorphism of $G_1$ preserving the labeling. Since we labeled $G_1$ in a distinguishing way, at first, so $f$ is the identity automorphism on $G_1$. Regarding to the labeling of copies of $G_2$ and since  $f$ preserves the labeling of the copies of $G_2$, so $f$ maps each copy of $G_2$ to itself. The map $f$ is the identity automorphism on each copy of $G_2$, because each copy of $G_2$ was labeled in a distinguishing way. Therefore $f$ is the identity automorphism on $G_1\star G_2$.\qed

 The following corollary is an immediate consequence of Theorem \ref{thm}.  
 \begin{corollary}
  Let $G_1$ and $G_2$ be two connected graphs of orders $n_1$ and $n_2$, respectively, such that $n_1 > 1$.  If $D(G_1)=1$, then $D(G_1\star G_2)\leqslant D(G_2)$.
 \end{corollary}
 \proof It is sufficient to note that if $D(G_1)=1$, then the value of $M$ in Theorem \ref{thm} is zero.\qed
 
 \medskip
 We end the paper by presenting  an upper bound for the distinguishing index of the neighbourhood corona  of two graphs: 
 \begin{theorem}
  Let $G_1$ and $G_2$ be two connected graphs of orders $n_1$ and $n_2$, respectively, such that $n_1 > 1$. Then $D'(G_1\star G_2)\leqslant max \{D'(G_1), D'(G_2)\}$.
 \end{theorem}
 \proof
 We  define an edge distinguishing labeling of $G_1\star G_2$ with  $max \{D'(G_1), D'(G_2)\}$ labels. To obtain such labeling we first label the edge set of $G_1$ and $G_2$ in a distinguishing way with $D'(G_1)$ and $D'(G_2)$ labels, respectively. For the labeling of the edges between each copy of $G_2$ and $G_1$ we use of the labeling of the edge set of $G_1$ as follows:
 
 Let $N_{G_1}(v_k)=\{v_{k1},\ldots , v_{1\vert N_{G_1}(v_k) \vert}\}$, where $1\leqslant k \leqslant n_1$. By the notations of the vertices of $G_1$ and the copies of $G_2$, we assign the all edges $v_{kj_k}u_i^k$, $1\leqslant i \leqslant n_2$, the label of the edge $v_{kj_k}v_k$ in the distinguishing labeling of the edge set of $G_1$, where $1\leqslant k \leqslant n_1$ and $1\leqslant j_k\leqslant \vert N_{G_1}(v_k) \vert$.  This labeling is a distinguishing edge labeling of $G_1\star G_2$, because if $f$ is an automorphism of $G_1\star G_2$ preserving the labeling, then by Corollary \ref{cor}, the restriction of  $f$ to $G_1$ is an automorphism of $G_1$ preserving the labeling. Since we labeled $G_1$ in a distinguishing way, at first, so $f$ is the identity automorphism on $G_1$. Regarding to the labeling of the edges between the copies of $G_2$ and $G_1$ and by Corollary \ref{cor} we conclude that  $f$ maps each copy of $G_2$ to itself. Since we labeled each copy of $G_2$ in a distinguishing way, at first, so the map $f$ is the identity automorphism on each copy of $G_2$, and so $f$ is the identity automorphism on $G_1\star G_2$.\qed
 

\begin{thebibliography}{99}

 
\bibitem{Albert} M.O. Albertson and K.L. Collins, {\it Symmetry breaking in graphs}, Electron. J. Combin. 3 (1996) \#R18.

 	\bibitem{soltani2} S. Alikhani and S. Soltani, {\it Distinguishing number and distinguishing index of certain graphs}, submitted. Available at \texttt{http://arxiv.org/abs/1602.03302}.
 
 
\bibitem{Frucht} R. Frucht and F. Harary, {\it On the corona two graphs}, Aequationes Math.  4  (1970) 322-325.


\bibitem{Harary} F. Harary, {\it Graph Theory}, Addition-Wesley Publishing Co., Reading, MA/Menlo Park, CA/London, 1969.

\bibitem{Gopalapillai}  I. Gopalapillai, {\it The spectrum of neighborhood corona of graphs}, Kragujevac Journal of Mathematics. 35 (2011) 493-500.


\bibitem{Kali1} R. Kalinowski and M. Pilsniak, {\it Distinguishing graphs by edge colourings}, European J. Combin. 45 (2015) 124-131.

	\bibitem{Klavzar} S. Kla\v{v}zar and X. Zhu, {\it Cartesian powers of graphs can be distinguished by two labels}, European J. Combin. 28 (2007) 303-310. 
	
	\bibitem{linear} X. Liu and S. Zhou, {\it Spectra of the neighbourhood corona of two graphs},  Linear  Multilinear Alg.  62, 9 (2014) 1205--1219.    


\bibitem{fish} F. Michael and I. Garth, {\it Distinguishing colorings of Cartesian products of complete graphs}, Discrete Math., 308 (11), (2008) 2240-2246. 

\bibitem{Sampathkumar} E. Sampathkumar, H. B. Walikar, {\it On the splitting graph of a graph}, Karnatak Univ. J. Sci. 35/36 (1980-1981), 13-16.
	
	 \end{thebibliography}
\end{document}